\documentclass[11pt,a4paper]{article}
\usepackage[margin=25mm]{geometry}
\usepackage[T1]{fontenc}
\usepackage{lmodern}
\usepackage{amsmath,amssymb,amsthm,mathtools}
\usepackage{booktabs,array}
\usepackage{xcolor}
\usepackage{float}
\usepackage[colorlinks=true,linkcolor=blue!45!black,citecolor=blue!45!black,urlcolor=blue!45!black]{hyperref}
\usepackage{todonotes}
 
\newtheorem{theorem}{Theorem}[section]
\newtheorem{proposition}[theorem]{Proposition}
\newtheorem{lemma}[theorem]{Lemma}

\theoremstyle{definition}

\newtheorem{remark}[theorem]{Remark}
\newtheorem*{theoremA}{Theorem A}
\newtheorem*{theoremB}{Theorem B}
\newtheorem*{corollaryC}{Corollary C}
\newcommand{\D}{\mathsf D}
\newcommand{\DD}{\mathsf D^{*}}
\newcommand{\F}{\mathbb F}
\newcommand{\Z}{\mathbb Z}

\newcommand{\GL}{\operatorname{GL}}

\newcommand{\ind}{\mathbf 1}
\newcommand{\It}{\mathcal I}
\numberwithin{equation}{section}
\title{Integral zero-sum incidence identities\\and Davenport constants in rank three}

\author{Claudiu G.~Pop\thanks{Mathematical Institute, Master ALGANT, Universit\'e de Bordeaux,
351 cours de la Lib\'eration,
33405 Talence, France.
E-mail: \texttt{popclaudiu38@gmail.com}} \and George C. Rad \and  George C.~\c{T}urca\c{s}\thanks{Faculty of Mathematics and Computer Science,
Babe\c{s}-Bolyai University, Cluj-Napoca, Romania.\\
E-mails: \texttt{george.ciprian.rad@stud.ubbcluj.ro},  \texttt{george.turcas@ubbcluj.ro}} }
\date{\today}
 
\begin{document}
\maketitle
\begin{abstract}
We study the Davenport constant of finite abelian groups through signed identities among zero-sum
subsequences. An integral incidence identity of degree $t$ for a sequence
assigns integer weights to its nonempty zero-sum subsequences so that the
weights add up to one, while those of the subsequences containing any
given nonempty set of at most $t$ terms add up to zero. Over a finite
abelian $p$-group $H$, the zero-sum subsequences of a sequence of any
length have an indicator polynomial of degree at most $\D(H)-1$, and this
gives a transfer principle: if, for each $t\ge0$, every sequence of $L(t)$
terms over $Q$ admits such an identity of degree $t$, then
$\D(Q\oplus H)\le L(\D(H)-1)$ for every finite abelian $p$-group $H$.
For $m\in\{2,3,4\}$ we construct integral identities at length
$L(t)=m(t+3)-2$ over $C_m^3$, and prove that these thresholds are sharp
over fields of characteristic zero or characteristic not dividing $m$. For a prime $p\nmid m$ and $0\le a\le b\le c$, we prove that
\[
\D(C_{mp^a}\oplus C_{mp^b}\oplus C_{mp^c})
=m(p^a+p^b+p^c)-2.
\]
In particular $\D(C_{15}^3)=43$ and $\D(C_{20}^3)=58$, which are the two smallest
cases of the conjecture $\D(C_n^3)=3n-2$ left open in the literature. The proofs of the base cases for $m \in \{3,4\}$ are computer-assisted. All of the results in this paper have been formalised in Lean 4. 
\end{abstract}

\section{Introduction}

Let $G$ be a finite abelian group, written additively. A \emph{sequence} over $G$ is a finite list of elements of $G$ in which repetitions are allowed, and the \emph{Davenport constant} $\D(G)$ is the least integer $L$ such that every
sequence of $L$ terms over $G$ has a nonempty subsequence with sum zero. If
\[
G\simeq C_{n_1}\oplus\cdots\oplus C_{n_r},
\qquad 1<n_1\mid\cdots\mid n_r,
\]
then $n_i-1$ copies of each basis element form a zero-sum free sequence, so
\[
\D(G)\ \ge\ \DD(G):=1+\sum_{i=1}^r(n_i-1).
\]
Olson proved that equality holds for $p$-groups and for groups of rank at
most two~\cite{OlsonI,OlsonII}.
Equality fails for infinitely many groups of rank four
\cite{GeroldingerSchneider,GaoGeroldinger}. In rank three, it is conjectured that
\begin{equation}\label{eq:conj}
\D(C_{n_1}\oplus C_{n_2}\oplus C_{n_3})=n_1+n_2+n_3-2
\qquad (1<n_1\mid n_2\mid n_3).
\end{equation}
For background on the conjecture, see
\cite{Gao2000,GaoGeroldinger,Schmid,GGZ}.
In this paper, we develop a method for attacking particular cases of this conjecture.
We write the group as $Q\oplus H$ with $H$ a $p$-group, and read a
sequence over it as a sequence over the \emph{quotient} $Q$ whose terms
carry labels in the \emph{kernel} $H$. We first study identities
among the zero-sum subsequences of sequences over a small quotient, and then use those
identities for arbitrary labels in a $p$-group kernel.

Our method allows us to prove the conjecture \eqref{eq:conj} for
infinitely many mixed-prime groups of rank three with a
quotient $C_3^3$ or $C_4^3$, and the values it produces include the two
smallest cases of~\eqref{eq:conj} for $C_n^3$ that the classical theory
leaves open, $n=15$ and $n=20$.

\paragraph{The inductive method.} The standard way to bound $\D(G)$ from above when $G$ is not a
$p$-group reduces the question to two smaller groups, a subgroup $H$ and
the quotient $Q=G/H$. It goes back to the earliest work on the subject
\cite{OlsonII,vanEmdeBoasII}, and is known in the zero-sum literature as
the \emph{inductive method} \cite[Section~3]{Schmid},
\cite{BhowmikEtAl,GirardSchmid}.
A zero sum in a quotient need not be a zero sum in the original group.
Project a sequence over $G$ to $Q$, and call a subsequence a
\emph{block} if its projection has sum zero.
Each block has a sum in $H$; if we can find
$\D(H)$ pairwise disjoint blocks, some nonempty collection of their
sums adds to zero in $H$, and the union of those blocks is a zero sum
in $G$ \cite[Proposition~2.6]{DelormeOrdazQuiroz}.
 
Short blocks make this argument efficient. Write $\eta(Q)$ for the least
length such that every sequence over $Q$ of that length contains a
nonempty zero-sum subsequence of length at most $\exp(Q)$.
Extracting short blocks successively gives, in particular,
\begin{equation}\label{eq:classical}
\D(Q\oplus H)\le \eta(Q)+(\D(H)-1)\exp(Q).
\end{equation}
The same estimate holds for $G$ in the extension above.
 
The last block need not be short, and this improves the count. With $\D(H)\ge2$, extract $\D(H)-1$ short blocks and then one
arbitrary zero-sum block. This requires at most
\[
\max\{\eta(Q)+(\D(H)-2)\exp(Q),\ \D(Q)+(\D(H)-1)\exp(Q)\}
\]
terms. The first expression ensures that the last short block can be
extracted; the second ensures that a final zero sum remains.
The maximum is the second expression precisely when
\[
\eta(Q)\le\D(Q)+\exp(Q).
\]
For $Q=C_p^2$, the values $\eta(Q)=3p-2$ and $\D(Q)=2p-1$ satisfy this
inequality \cite{OlsonII,GaoGeroldinger}; this is the count that makes
the inductive method exact in rank two. For $Q=C_3^3$ and $Q=C_4^3$, however, the
comparison is $17>7+3$ and $22>10+4$, respectively
\cite{ShortII,Harborth}. For the groups $C_m^3\oplus H$ studied below,
the second expression is the basis lower bound, so for $m=3,4$ the count
overshoots it. Nor is the loss due to the insistence on short blocks: ten
terms over $C_3^3$ need not contain two disjoint blocks of any size
(Appendix~\ref{app:computations}), although $\D(C_3^3\oplus C_2)=10$.
 
Once disjoint blocks have been chosen, this argument uses only their sums
in $H$. It forgets the other zero sums of the projection and the ways
in which they overlap. Our approach keeps that incidence information
and combines the zero sums with signed weights.

\paragraph{Incidence identities.}
We regard sequences as indexed, so equal group elements at different
positions are considered different terms, since the kernel labels attached to
those positions may differ. For $U=(u_1,\dots,u_N)$ over $G$, write
\[
\mathcal Z(U)=
\{\,\varnothing\ne E\subseteq[N]:\textstyle\sum_{j\in E}u_j=0\,\}
\]
for the family of its \emph{zero sets}, and let $\ind_E{}\in\{0,1\}^N$ be the incidence vector
of $E$.
 
Let $t$ be a non-negative integer, and let $R$ be a commutative ring with identity.
An \emph{incidence identity of degree $t$ over $R$}, or
\emph{$\It_t$-identity}, is a family $(c_E)_{E\in\mathcal Z(U)}$
of elements of $R$ satisfying
\begin{equation}\label{eq:moments}
\sum_E c_E=1,\qquad
\sum_{E\supseteq J}c_E=0
\quad\text{for every }J\subseteq[N]\text{ with }1\le |J|\le t.
\end{equation}
The sum indexed by $J$ is the \emph{moment} at $J$.
Thus the weights add to one, but the weights through any prescribed
nonempty set of at most $t$ positions cancel. We say the identity is
\emph{integral} when $R=\Z$.
 
For example, if $B$ and $C$ are disjoint nonempty zero
sets, give weights $1,1,-1$ to $B,C,B\cup C$, respectively.
The total weight is one, and the total through each individual position
is zero. These are exactly the conditions for degree one.
 
The polynomial meaning is the reason for requiring all these
cancellations. On Boolean inputs, the monomial
$\prod_{j\in J}z_j$ records whether a set contains every position of $J$.
Consequently~\eqref{eq:moments} is equivalent to
\[
\sum_E c_E P(\ind_E)=P(0)
\quad\text{for every polynomial $P$ of degree at most $t$}.
\]
We prove this in Lemma~\ref{lem:duality}. The empty subsequence receives
no weight; nevertheless, these signed evaluations reproduce its value
for every polynomial of the indicated degree.
 
For an abelian $p$-group $H$, the test for a selected set of kernel
labels to have sum zero is such a polynomial over $\F_p$: its degree
is at most $\D(H)-1$, independently of the number of labels
(Lemma~\ref{lem:indicator}). An identity of this degree for the
sequence over the quotient can therefore be evaluated at the kernel's
zero-sum test, and this gives the transfer principle.

\begin{theoremA}[Transfer principle]
Let $Q$ be a finite abelian group, let $p$ be a prime, and let
$L(0),L(1),L(2),\dots$ be integers such that every sequence of $L(t)$ terms
over $Q$ has an $\It_t$-identity over $\F_p$. Then
\[
\D(Q\oplus H)\ \le\ L\bigl(\D(H)-1\bigr)
\qquad\text{for every finite abelian $p$-group } H.
\]
If the identities are integral, the conclusion holds for every prime $p$.
\end{theoremA}

To see the mechanism, let $P$ be the zero-sum indicator for the kernel
labels and evaluate the quotient identity at $P$. The result is
\[
1=P(0)=\sum_E c_E P(\ind_E).
\]
Some summand is nonzero, so some nonempty zero set in the quotient also
has sum zero in the kernel. For a fixed degree, the weights depend only
on the quotient sequence; the same identity works for every choice of
kernel labels. Integer weights can be reduced modulo any prime, which
makes a single certificate useful for many different kernels.
 
Over a field $K$, the alternative to an $\It_t$-identity is a
\emph{separator of degree at most $t$}: a polynomial $P$ of degree at most
$t$ such that $P(0)=1$ and $P(\ind_E)=0$ for every $E\in\mathcal Z(U)$.
Lemma~\ref{lem:duality} shows that exactly one of these alternatives
occurs. We write $\kappa_t(G;K)$ for the least $N$ such that no
sequence of $N$ terms over $G$ has a separator of degree at most $t$.
In degree zero the only possible separator is the constant polynomial
$1$, which works exactly when there is no nonempty zero sum. Thus
$\kappa_0(G;K)=\D(G)$.
 
\paragraph{The quotient hierarchy.}
For $Q=C_m^3$ and
$H=C_{p^a}\oplus C_{p^b}\oplus C_{p^c}$ with $p\nmid m$, Olson's theorem
and the Chinese remainder theorem give
\[
\DD(C_m^3\oplus H)=m(\D(H)-1)+3m-2.
\]
This singles out the target length
$L_m(t)=mt+3m-2=m(t+3)-2$: a quotient identity at this length has
exactly the degree and length needed to reach the basis lower bound
after transfer.

\begin{theoremB}
Let $m\in\{2,3,4\}$ and let $t\ge0$ be an integer. Every indexed sequence of
$m(t+3)-2$ terms over $C_m^3$ has an integral $\It_t$-identity.
Consequently $\kappa_t(C_m^3;K)\le m(t+3)-2$ for every field $K$, with
equality whenever the characteristic of $K$ is zero or does not divide~$m$.
\end{theoremB}

The $m=2$ case gives a small geometric model of the argument: a signed
identity on the Fano plane with one point removed. For $m=3$ and $m=4$,
finite computations supply the base cases, and the three constructions
of Lemma~\ref{lem:operations} propagate them to all degrees.
 
Combining Theorems A and B gives
$\D(C_m^3\oplus H)\le m(\D(H)-1)+3m-2$
for every finite abelian $p$-group $H$. When $H$ has rank at most
three and $p\nmid m$, the right-hand side is the basis lower bound, giving the
following consequence.

\begin{corollaryC}
Let $m\in\{2,3,4\}$, let $p$ be a prime with $p\nmid m$, and let
$0\le a\le b\le c$ be integers. Then
\begin{equation}\label{eq:family}
\D(C_{mp^a}\oplus C_{mp^b}\oplus C_{mp^c})
=m(p^a+p^b+p^c)-2.
\end{equation}
In particular $\D(C_{3p^a}^3)=9p^a-2$ for $p\ne3$ and
$\D(C_{4p^a}^3)=12p^a-2$ for odd~$p$.
\end{corollaryC}

For example,
$\D(C_{15}^{3})=43$. As $C_{15}^3\simeq C_3^3\oplus C_5^3$ and $\D(C_5^3)=13$, the classical estimate~\eqref{eq:classical} only gives the upper bound $53$. 
 
Trivial kernel factors are allowed, with $C_{p^0}=C_1$.
When $p\mid m$, the rank-three groups displayed in~\eqref{eq:family}
are themselves $p$-groups, and their Davenport constants are already
given by Olson's theorem.

\paragraph{Earlier work and our contribution.}
Corollary C includes classical families: the doubled $p$-groups and
tripled $2$-groups of van Emde Boas and Kruyswijk
\cite{vanEmdeBoasII,vanEmdeBoasIII}, as well as special cases covered by
later rank-three results \cite{BhowmikSchlagePuchta,BhowmikEtAl,Schmid,GirardSchmid}.
For $C_n^3$ earlier work covers, besides the prime powers, the
exponents $n=2p^a$ and $n=3\cdot2^a$ \cite{Schmid,Girard}, but not
$n=3p^a$ or $n=4p^a$ with $p\ge5$.\footnote{A remark preceding
\cite[Theorem~6.6]{GeroldingerYang} lists $n=3p^a$ among the known cases,
citing \cite[Section~4]{Schmid}. The tripling result recalled there
requires a group of even order and gives $n=3\cdot2^a$; we know of no
proof for $p\ge5$.}
 
The group-algebra input goes back to Olson, and polynomial methods have
a substantial history in zero-sum theory \cite{OlsonI,Alon}.
Finite computational input to the inductive method also occurs in the work of
Bhowmik, Halupczok and Schlage-Puchta \cite{BhowmikEtAl}.
Related signed zero-block and moment methods occur in Volkmann's work
on nonabelian groups of exponent three
\cite{Volkmann2026a,Volkmann2026b,Volkmann2026c}.
Here we establish complete incidence identities over the integers,
throughout the degree hierarchy for three small quotients.
The transfer principle explains how those identities
control arbitrary labels in an abelian $p$-group kernel, while
integrality supplies a common certificate in every characteristic.
 
\paragraph{Scope.}
The bounded-degree indicator in Theorem A uses a $p$-group kernel.
Remark~\ref{rem:limit} explains why a degree bound independent of the
number of labels fails for other kernels in characteristic $p$.
If $n$ has at least three distinct prime divisors, separating off one
primary component of $C_n^3$ leaves a quotient $C_m^3$ with at least two
prime divisors. An incidence hierarchy at length $m(t+3)-2$ for such a quotient would already
include the Davenport conjecture for that quotient in degree zero.
Deleting short zero sets, as in Section~\ref{sec:bases}, reduces a proposed hierarchy to finitely many
base degrees for each fixed $m$; establishing those base cases is still
a mathematical task.
 
\paragraph{Organisation.}
Section~\ref{sec:identities} develops the identities and their three
constructions. Section~\ref{sec:transfer} proves the transfer principle,
and Section~\ref{sec:sharp} gives the lower bound.
Section~\ref{sec:bases} establishes the quotient hierarchies and
Section~\ref{sec:davenport} is dedicated to the proof of the announced Davenport values.
Section~\ref{sec:data} records the computational dependencies and
the formal verification, with the
coverage arguments and certificate specifications in
Appendix~\ref{app:computations}.

\paragraph{AI disclosure.} Large language models (OpenAI's Codex and ChatGPT, Anthropic's Claude) were used in an essential way for this work: in the exploration that led to the incidence identities, which play a key role in the main results of this paper, in coding the enumeration programs, in writing and revising the exposition and the proofs, and in a substantial part of the Lean formalisation.  The authors made a significant effort to improve the exposition of the AI-produced mathematics and to give proper credit to prior work that uses similar ideas.  Every argument was check independently by at least one author and the author take full responsibility for the paper.

\section{Incidence identities and three constructions}\label{sec:identities}
 
For a sequence $U=(u_1,\dots,u_N)$ write
$\sigma(U)=\sum_{j=1}^N u_j$; the notation $x^k$ denotes $k$ occurrences
of $x$. A subsequence is specified by its set of positions.

\begin{lemma}\label{lem:duality}
Let $(c_E)$ be a family of elements of $R$ indexed by $\mathcal Z(U)$.
\begin{enumerate}
\item $(c_E)$ is an $\It_t$-identity over $R$ if and only if
\begin{equation}\label{eq:evaluation}
\sum_{E\in\mathcal Z(U)}c_E\,P(\ind_E)=P(0)
\end{equation}
for every polynomial $P\in R[z_1,\dots,z_N]$ of degree at most $t$.
\item An $\It_t$-identity over $\Z$ is an $\It_t$-identity over every
commutative ring, by reduction of the coefficients.
\item If $R=K$ is a field, $U$ has an $\It_t$-identity over $K$ if and only
if $U$ has no separator of degree at most $t$ over $K$.
\end{enumerate}
\end{lemma}
\begin{proof}
Both sides of~\eqref{eq:evaluation} are linear in $P$. On Boolean inputs,
every polynomial agrees with a multilinear polynomial of no greater
degree, so it is enough to test the monomials
$z_J=\prod_{j\in J}z_j$ with $|J|\le t$, including $z_\varnothing=1$.
For these monomials the equality is
$\sum_{E\supseteq J}c_E=[J=\varnothing]$, exactly the moment conditions.
 
Part (2) is clear, because we always have a ring morphism $\mathbb Z\to R$. 
 
For (3), consider the vectors
$m_E=(\ind_E{}^J)_{|J|\le t}\in K^{\binom{N}{\le t}}$ of monomial values, for
$E\in\mathcal Z(U)$, and the vector $m_\varnothing=(1,0,\dots,0)$ of the
empty set. Observe that an identity is an expression of $m_\varnothing$ as a linear
combination of the $m_E$, and a separator is a linear functional vanishing
on all $m_E$ and equal to $1$ on $m_\varnothing$.
Over a field, either $m_\varnothing$ belongs to the span of the vectors
$m_E$, or there is a linear functional vanishing on that span and taking
the value $1$ at $m_\varnothing$. The first alternative gives the identity;
the second gives the separator.
\end{proof}

An identity on a subsequence extends to the full sequence by assigning
zero weight to the other zero sets. Conversely, a separator for the
full sequence restricts to one for a subsequence by setting the omitted
coordinates to zero. Thus existence of an identity, and absence of a
separator, are inherited by longer sequences. The invariant $\kappa_t$
is finite: $(t+1)\D(G)$ terms contain $t+1$ disjoint zero sets, and
Lemma~\ref{lem:operations}(1) supplies an identity.
 
The identities that Section~\ref{sec:bases} obtains by computer are
stored with integer weights. Reducing their
coefficients modulo $p$ gives the identities needed for Theorem A, so one
exact certificate check serves every characteristic. Solvability
over each field separately would give neither the integral assertion of
Theorem~B nor such a common certificate.

We use three ways to build identities. Disjoint zero sets give
inclusion--exclusion. Adjoining a disjoint zero set turns a finite
difference into an identity of one higher degree. Taking complements
inside a zero set cancels the highest-degree terms in a different way.
 
\begin{lemma}[Three constructions]\label{lem:operations}
Let $U$ be an indexed sequence and $t\ge0$. The following produce integral
$\It_t$-identities for $U$ from integral identities of lower degree.
\begin{enumerate}
\item If $U$ has $t+1$ pairwise disjoint zero sets, then $U$ has an
$\It_t$-identity.
\item Let $t\ge1$. If $B\in\mathcal Z(U)$ and a subsequence of $U$ disjoint
from $B$ has an $\It_{t-1}$-identity, then $U$ has an $\It_t$-identity.
\item Let $t\ge1$. If $Z\in\mathcal Z(U)$ and a proper subsequence
$T\subsetneq Z$ has an $\It_{t-1}$-identity, then $U$ has an $\It_t$-identity.
\end{enumerate}
\end{lemma}

\begin{proof}
For (1), let $B_1,\dots,B_{t+1}$ be disjoint nonempty zero sets.
Give the union of every nonempty collection of $k$ blocks the weight
$(-1)^{k+1}$, and give every other zero set weight zero.
Distinct collections have distinct unions. If a nonempty set $J$ of
at most $t$ positions is contained in such a union, it misses some block
$B_i$. The collections whose union contains $J$ are paired by adding or
removing $B_i$, with opposite weights. Thus every required nonconstant
moment vanishes. The total weight is
$\sum_{k=1}^{t+1}(-1)^{k+1}\binom{t+1}{k}=1$.
 
For (2), let $(c_E)$ be the identity of degree $t-1$ on a subsequence
$U'$ disjoint from $B$. Define new weights by
\[
d_E=c_E,\qquad d_{E\cup B}=-c_E
\quad(E\in\mathcal Z(U')),\qquad d_B=1,
\]
and set the remaining weights to zero. These sets are nonempty zero
sets and are pairwise distinct. Let $P$ be multilinear of degree at
most $t$. Restricting to the coordinates outside $B$, the polynomial
\[
P(z+\ind_B)-P(z)
\]
has degree at most $t-1$. Indeed, terms involving no coordinate of $B$
cancel, and each remaining term has spent at least one degree on a
coordinate of $B$. Applying the old identity gives
\[
\sum_{E\in\mathcal Z(U')}c_E
  \bigl(P(\ind_{E\cup B})-P(\ind_E)\bigr)
=P(\ind_B)-P(0).
\]
Rearranging yields $\sum_F d_F P(\ind_F)=P(0)$.
 
For (3), let $(c_E)$ be the identity on $T\subsetneq Z$, and put
$s=(-1)^t$. Give weight $c_E$ to $E$, weight $-s c_E$ to $Z\setminus E$,
and weight $s$ to $Z$, with zero weights elsewhere.
The complements are nonempty zero sets. The assigned sets are distinct:
each complement contains $Z\setminus T$, whereas each $E$ lies in $T$;
and $Z$ is neither an $E$ nor the complement of a nonempty $E$.
Restrict a multilinear polynomial $P$ of degree at most $t$ to the
coordinates of $Z$. Then
\[
R_P(z)=P(z)-sP(\ind_Z-z)
\]
has degree at most $t-1$: complementation multiplies each degree-$t$
term by $(-1)^t$, so the degree-$t$ parts cancel.
Applying the identity on $T$ to $R_P$ and rearranging gives
\[
\sum_E c_E P(\ind_E)
-s\sum_E c_E P(\ind_{Z\setminus E})
+sP(\ind_Z)=P(0).
\]
All three constructions use only addition and subtraction of coefficients,
so they preserve integrality.
\end{proof}

\section{The transfer principle}\label{sec:transfer}
 
For a cyclic kernel $H=C_p$, identified with the additive group of
$\F_p$, the required indicator is familiar. Fermat's theorem gives
\[
P(z)=1-\Bigl(\sum_j v_j z_j\Bigr)^{p-1}.
\]
At a Boolean vector $z=\ind_E$, this is $1$ precisely when
$\sum_{j\in E}v_j=0$, and it has degree at most $p-1$, however many labels
there are. Multiplying one such factor for each coordinate gives the
same phenomenon for $C_p^r$. The next lemma extends it to arbitrary
abelian $p$-groups, including cyclic factors of order $p^a$.
 
\begin{lemma}\label{lem:indicator}
Let $H=\bigoplus_{i=1}^r C_{p^{a_i}}$ and $M=\sum_i(p^{a_i}-1)$. For
arbitrary $v_1,\ldots,v_N\in H$, the Boolean function
\[
f(z)=\ind_{\{\sum_j z_jv_j=0\}}\in\F_p ,\qquad z\in\{0,1\}^N,
\]
is represented by a multilinear polynomial over $\F_p$ of degree at most $M$.
\end{lemma}

\begin{proof}
Write $X^h$ for the basis element of $\F_p[H]$ indexed by $h\in H$,
so that $X^aX^b=X^{a+b}$. Choose generators $e_i$ of the cyclic factors,
with $e_i$ of order $p^{a_i}$, and put $Y_i=X^{e_i}-1$.
Then
\[
\F_p[H]\simeq
\F_p[Y_1,\ldots,Y_r]/(Y_1^{p^{a_1}},\ldots,Y_r^{p^{a_r}}),
\]
and the augmentation ideal $I=(Y_1,\ldots,Y_r)$ satisfies $I^{M+1}=0$, since
a monomial of total degree $M+1$ in the $Y_i$ has some exponent at least
$p^{a_i}$. Now we consider the multilinear polynomial
\[
P(z)=[X^0]\prod_{j=1}^N\bigl(1+z_j(X^{v_j}-1)\bigr),
\]
where $[X^0]$ extracts the identity coefficient in the group basis.
Multiplication records the selected sum: setting $z_j=0$ omits $v_j$,
and setting $z_j=1$ includes it. For a Boolean vector
$\ind_E$ the product collapses to $X^{\sum_{j\in E}v_j}$, so $P(\ind_E)=f(\ind_E)$. The coefficient of $z_J$ in $P$ is
$[X^0]\prod_{j\in J}(X^{v_j}-1)$. Before taking $[X^0]$, this is a
product of $|J|$ elements of the augmentation ideal, and therefore lies in
$I^{|J|}$. Once $|J|>M$, that product, and hence the corresponding scalar
coefficient of $P$, is zero. This is why the degree is bounded by $M$
rather than by the number $N$ of labels. Thus the degree of our polynomial $P$ is at most $M$.
\end{proof}

\paragraph{Olson's theorem.}
The same nilpotence gives $\D(H)=M+1$. If $M+1$ terms were zero-sum
free, the product $\prod_j(1-X^{v_j})$ would have identity coefficient
$1$, contributed only by the empty selection. But this product lies
in $I^{M+1}=0$. Hence $\D(H)\le M+1$.
Conversely, $p^{a_i}-1$ copies of each $e_i$ form a zero-sum free
sequence of length $M$, giving equality.
For this basis sequence, only the empty subset has sum zero. Its unique
multilinear indicator is $\prod_{j=1}^M(1-z_j)$, of degree $M$.
Thus the indicator degree bound is best possible in general.

\begin{proof}[Proof of Theorem~A]
Put $M=\D(H)-1$ and take $L(M)$ terms $(u_j,v_j)$ in $Q\oplus H$.
The quotient sequence $(u_j)$ has an $\It_M$-identity $(c_E)$ over
$\F_p$. Let $P$ be the indicator polynomial for the kernel labels
$(v_j)$ supplied by Lemma~\ref{lem:indicator}. Since $\deg P\le M$,
\[
1=P(0)=
\sum_{E\in\mathcal Z((u_j))}c_E P(\ind_E).
\]
Some summand is nonzero. For its nonempty set $E$,
$\sum_{j\in E}u_j=0$ by definition, while $P(\ind_E)=1$ gives
$\sum_{j\in E}v_j=0$. Thus the original sequence has a nonempty zero sum,
and $\D(Q\oplus H)\le L(M)$.
If the identities are integral, reducing their coefficients modulo $p$
makes the same argument available for every prime.
\end{proof}

\begin{remark}[The elementary block bound]\label{rem:classical}
A sequence of length $\eta(Q)+t\exp(Q)$ contains $t+1$ disjoint short
zero sets, obtained by repeated extraction. Their unions, with
inclusion--exclusion weights, give an $\It_t$-identity by
Lemma~\ref{lem:operations}(1). For a $p$-group kernel, Theorem A therefore
recovers~\eqref{eq:classical}, and in the same way every bound
obtained from $\D(H)$ disjoint blocks, whatever their sizes. Incidence identities allow further
cancellations among zero sets that do not arise from a disjoint block
decomposition. The six-term Fano example below already requires such
overlap.
\end{remark}

\begin{remark}[The kernel hypothesis]\label{rem:limit}
The degree bound in Lemma~\ref{lem:indicator} cannot hold uniformly in
the number of labels when $H$ is not a $p$-group. Choose an element $g$
of prime order $q\ne p$ and take all $N=p^s$ labels equal to $g$.
The coefficient of $z_1\cdots z_N$ in their zero-sum indicator is
\[
[X^0](X^g-1)^N
=[X^0](X^{Ng}-1)=-1\ne0,
\]
where the middle equality is Frobenius and $Ng\ne0$ because $q\nmid N$.
A function on $\{0,1\}^N$ has a unique multilinear representative over
$\F_p$, so the indicator has degree exactly $N$. Since $s$ is arbitrary,
there is no bound depending only on $H$.
 
For example, when $g$ has order two and $p$ is odd, the indicator is
\[
\frac12\left(1+\prod_{j=1}^N(1-2z_j)\right),
\]
the test that the number of selected labels is even.
Thus the bounded-degree indicator argument cannot be applied unchanged
to such kernels.
\end{remark}

\section{Sharpness}\label{sec:sharp}
 
\begin{proposition}\label{prop:sharp}
Let $m\ge2$ and $t\ge0$, and let $K$ be a field of characteristic zero or of
characteristic not dividing $m$. Then
\[
\kappa_t(C_m^3;K)\ge m(t+3)-2 .
\]
\end{proposition}

\begin{proof}
We choose a sequence whose zero sets are especially simple:
\[
U=e_1^{\,m-1}e_2^{\,m-1}e_3^{\,m(t+1)-1},
\]
of length $m(t+3)-3$.
A zero set contains no first- or second-axis term, because there are
fewer than $m$ of each. It therefore consists of third-axis positions,
and its size is $mk$ for some $1\le k\le t$.
 
We look for a separator among the linear combinations of the
elementary symmetric polynomials
$s_0,\dots,s_t$ in those third-axis coordinates. On a set of size $l$,
their values are $s_j(\ind_E)=\binom{l}{j}$. It suffices to show that the
integer matrix
\[
B=\left(\binom{mk}{j}\right)_{0\le k,j\le t}
\]
is invertible over $K$.
 
To compute its determinant, let $a_{j,i}$ count the $j$-element subsets
of a fixed union of $i$ disjoint $m$-element blocks that meet every block.
Put $a_{0,0}=1$ and $a_{j,0}=0$ for $j>0$.
A $j$-element subset of a union of $k$ blocks first chooses the $i$ blocks it meets, giving
\[
\binom{mk}{j}=\sum_{i=0}^{k}\binom{k}{i}a_{j,i}.
\]
Thus $B=ST$, where
$S=(\binom{k}{i})_{0\le k,i\le t}$ is lower triangular with diagonal
entries one, and $T=(a_{j,i})_{0\le i,j\le t}$ is upper triangular.
Moreover, $a_{j,j}=m^j$: a $j$-element subset meeting $j$ blocks chooses
one element from each. Consequently
\[
\det B=\prod_{j=0}^{t}m^j=m^{t(t+1)/2}.
\]
This identity was obtained over the integers without division.
The determinant has nonzero image in $K$ under the stated characteristic
hypothesis. We can therefore choose a linear combination of
$s_0,\dots,s_t$ taking the values $1,0,\dots,0$ at subset sizes
$0,m,\dots,tm$. It is a separator for $U$ of degree at most $t$.
Hence $\kappa_t(C_m^3;K)\ge m(t+3)-2$.
\end{proof}

\begin{remark}
    The restriction on the characteristic cannot be dropped, and
the smallest case shows the threshold moving: $\kappa_1(C_2^3;\F_2)=5$,
one less than the value $6$ that Theorem~B gives over every field of
characteristic different from two. Indeed, for five terms $u_1,\dots,u_5$
over $C_2^3=\F_2^3$ the linear map $\F_2^5\to\F_2^3$,
$z\mapsto\sum_j z_ju_j$, has a null space of dimension at least two, and its
nonzero elements are the incidence vectors of the zero sets. If $x$ and $y$
are two independent vectors of this null space, an affine separator $1+\ell$ with
$\ell$ linear would have to satisfy $\ell(x)=\ell(y)=\ell(x+y)=1$, which is
impossible; so no five-term sequence has an affine separator over $\F_2$.
On the other hand $e_1e_2e_3^2$ has the affine separator $1+z_3$, its only
zero set being $\{3,4\}$.
Thus $\kappa_1(C_2^3;\F_2)=5$.
The integral six-term identity supplied by Theorem B remains valid in
characteristic two; the least length is simply smaller there.
\end{remark}
 
\section{The quotients \texorpdfstring{$C_2^3$, $C_3^3$ and $C_4^3$}{C2 cubed, C3 cubed and C4 cubed}}\label{sec:bases}

We now prove that every sequence of
$L_m(t)=m(t+3)-2$ terms over $C_m^3$ has an integral $\It_t$-identity.
Call a sequence \emph{short-zero-free} if it has no nonempty zero set of
size at most the exponent. For $t\ge1$, if a sequence of length $L_m(t)$ has a zero
set $B$ with $|B|\le m$, deleting it leaves at least
\[
L_m(t)-m=L_m(t-1)
\]
terms. An identity of degree $t-1$ on the remainder therefore gives one
of degree $t$ by Lemma~\ref{lem:operations}(2).
We call this step the \emph{short-zero recurrence}.
Once $L_m(t)\ge\eta(C_m^3)$, it is automatic.
The work is concentrated in the degrees before that point.
 
Identities of degree one and two are tested on affine and on
quadratic polynomials, and we call them \emph{affine} and \emph{quadratic}
identities. The following table summarises where the finite work lies;
the computations behind the finite certificates are described in
Appendix~\ref{app:computations}, and sharpness is
Proposition~\ref{prop:sharp}.
 
\begin{center}\small
\begin{tabular}{llp{6.4cm}l}
\toprule
Quotient & Target length & Before the automatic recurrence & Automatic recurrence\\
\midrule
$C_2^3$ & $2t+4$ & $t=0$: $\D(C_2^3)=4$; $t=1$: the Fano configuration & $t\ge2$, using $\eta(C_2^3)\le8$\\
$C_3^3$ & $3t+7$ & $t=0$: $\D(C_3^3)=7$; $t=1$: finite affine certificates; $t=2,3$: the support facts and the constructions & $t\ge4$, using $\eta(C_3^3)\le17$\\
$C_4^3$ & $4t+10$ & $t=0$: $\D(C_4^3)=10$; $t=1,2$: finite certificates & $t\ge3$, using $\eta(C_4^3)\le22$\\
\bottomrule
\end{tabular}
\end{center}

\subsection{The case \texorpdfstring{$C_2^3$}{C2 cubed}} 
The Fano configuration shows the mechanism in miniature.
Degree zero is $\D(C_2^3)=4$. For degree one, consider six terms.
If a term is zero or repeated, there is a zero set of at most two
positions; deleting it leaves at least four terms and hence a second,
disjoint zero set. Lemma~\ref{lem:operations}(1) applies.
 
Otherwise the six terms are distinct nonzero vectors of $\F_2^3$:
the Fano plane with one point $x$ removed. Its lines are the triples
$\{u,v,u+v\}$, each of which has sum zero. Give weight $1$ to the four
lines not through $x$, and weight $-1$ to the complements, in the
seven-point plane, of the three lines through $x$.
The latter are zero quadruples and do not contain $x$.
Each remaining point lies on two of the selected lines and in two of
the selected quadruples. Its total weight is therefore $2-2=0$,
whereas the total weight of the seven sets is $4-3=1$.
 
This identity uses overlap in an essential way. Every zero set among
six distinct nonzero vectors has at least three terms. Two disjoint
zero sets would therefore exhaust the six points, but their total is
the omitted nonzero point $x$, since all seven Fano points sum to zero.
 
Finally, any eight terms contain zero or a repetition, so
$\eta(C_2^3)\le8$. The short-zero recurrence gives all degrees $t\ge2$
at length $2t+4$. These are integral identities, valid after reduction
in every characteristic; Proposition~\ref{prop:sharp} gives equality
in characteristic zero and in odd characteristic.

\subsection{The case \texorpdfstring{$C_3^3$}{C3 cubed}}

We need the ten-term incidence identity and two elementary facts about supports.
 
\begin{lemma}[Ternary base facts]\label{lem:ternary} The following assertions hold.
\begin{enumerate}
\item Every ten-term sequence over $C_3^3$ has an integral $\It_1$-identity.
\item The support of a short-zero-free sequence over $C_3^3$ has at most
eight points. Every such support with seven points has nonzero sum, and every
such support with eight points has sum zero.
\end{enumerate}
\end{lemma}

We start by proving the second assertion.

\begin{proof}[Proof of Lemma~\ref{lem:ternary}(2)]
Let $A$ be the support of a short-zero-free sequence and put
$C=A\cup\{0\}$. Then $C$ is a \emph{cap}, that is, it contains no three distinct collinear points
(this is the connection with affine caps studied in~\cite{EdelEtAl}).
Indeed, over $\F_3$, three distinct points
are collinear exactly when their sum is zero. A line through $0$ would
give an opposite pair in $A$, and any other line would give a zero
triple in $A$.
 
An affine plane contains at most four cap points. Otherwise choose five.
In each of the four directions, the three parallel lines would contain
$2,2,1$ of these points, accounting for two pairs in that direction.
This accounts for at most eight pairs, rather than $\binom52=10$.
 
Put $N=|C|$ and let $n_i$ be the number of affine planes meeting $C$ in
exactly $i$ points, for $0\le i\le4$.
There are $(3^3-1)/2=13$ parallel classes of planes, each with three
members, hence $39$ planes. A point belongs to $13$ planes, a pair to
$(3^2-1)/2=4$, and a noncollinear triple to one.
Counting these incidences gives
\[
\begin{aligned}
\sum_{i=0}^4 n_i&=39,&
\sum_{i=0}^4 i n_i&=13N,\\
\sum_{i=0}^4\binom{i}{2}n_i&=4\binom N2,&
\sum_{i=0}^4\binom{i}{3}n_i&=\binom N3.
\end{aligned}
\]
For $N=10$, these equations give $n_3+4n_0=-4$, a contradiction.
A larger cap would contain a ten-point cap, so $N\le9$ and $|A|\le8$.
 
Suppose $|A|=7$, so $N=8$. The same equations give
$n_2=6n_0-2$, hence some affine plane is disjoint from $C$.
Write it as $\ell=c$, where $\ell$ is a nonzero linear functional.
Since $0\in C$, one has $c\ne0$. The other two parallel planes contain
four points each, and therefore
\[
\ell\left(\sum_{x\in C}x\right)
=4\sum_{a\in\F_3\setminus\{c\}}a=-c\ne0.
\]
Thus $\sigma(A)=\sigma(C)\ne0$.
 
Finally suppose $|A|=8$, so $N=9$.
No affine plane is empty, because the other two parallel planes could
hold at most eight cap points. Thus $n_0=0$, and the incidence equations
give $n_2=6n_0=0$. For every nonzero linear functional $\ell$, the three
parallel slices have sizes $3,3,3$ or a permutation of $1,4,4$.
In either case their sizes are congruent modulo three, so
$\ell(\sigma(C))=0$. This holds for every $\ell$, hence
$\sigma(A)=\sigma(C)=0$.
\end{proof}
 
\begin{proof}[Proof of Lemma~\ref{lem:ternary}(1), computer-assisted]
If a ten-term sequence has a zero set of at most three positions, its
deletion leaves at least seven terms and hence a disjoint zero set,
since $\D(C_3^3)=7$. Lemma~\ref{lem:operations}(1) then gives the
required affine identity.
 
It remains to consider the short-zero-free sequences. Their
multiplicities are at most two, and their supports satisfy the cap
conditions used above. The complete enumeration and its symmetry
reductions are given in Appendix~\ref{app:computations}.
Each case is supplied either with two disjoint zero sets or with an
explicit integral affine identity. The certificate check verifies the
zero sums, total weight one, and zero total weight through each
individual position. The transports between normalised representatives
and the indexed sequences they represent are checked together with the
certificates.
\end{proof}

Call a set $A\subseteq C_3^3$ \emph{allowed} if it contains neither zero
nor a pair of opposite points, and no three distinct points of $A$
sum to zero. Write $2A$ for the sequence containing two copies of
each point of $A$.
 
\begin{lemma}\label{lem:longternary}
Every short-zero-free sequence of thirteen terms over $C_3^3$ has a zero set
of size at least eleven.
\end{lemma}
\begin{proof}
The plan is to write $\sigma(U)$ as a sum of at most two terms of $U$:
deleting them leaves a zero set of size at least eleven. The representation
will come from the maximality of an allowed support, and the only point to
watch is whether it uses a term more often than $U$ contains it.
 
We first observe that if $A$ is an inclusion-maximal allowed support, then
every $x\in C_3^3$ is a sum of at most two terms of $2A$. For $x=0$ take no
terms, for $x\in A$ one term, and if $-x\in A$ take two copies of $-x$.
Otherwise $-x$ is a nonzero point outside $A$ whose opposite $x$ is also
outside $A$, so adjoining $-x$ to $A$ can only violate the third condition:
there is a zero triple $-x+a+b=0$ with distinct $a,b\in A$, which represents
$x=a+b$.
 
Now let $U$ be short-zero-free of length thirteen. Its multiplicities are at
most two, so its support has seven or eight points, and there are three
cases.
 
\emph{Case 1: the support is an eight-set $B$.} Then $U=2B-b-c-d$ for three
distinct points $b,c,d$ of multiplicity one. By Lemma~\ref{lem:ternary},
$\sigma(B)=0$, so $\sigma(U)=-b-c-d$, and $B$ is maximal by the support
bound. A representation of $\sigma(U)$ by at most two terms of $2B$ is
unavailable in $U$, that is, it uses a term more often than $U$ contains it, only if it uses one of $b,c,d$ twice, say $2b=-b-c-d$;
this gives $c+d=0$, impossible in an allowed support.
 
\emph{Case 2: the support is a seven-set $A$ which is maximal.} Then
$U=2A-a_0$ for one point $a_0$ of multiplicity one, and the only unavailable
representation is $2a_0=\sigma(U)=2\sigma(A)-a_0$, which forces
$\sigma(A)=0$, contradicting Lemma~\ref{lem:ternary}.
 
\emph{Case 3: the support is a seven-set $A$ which is not maximal.} Again
$U=2A-a_0$, and $A$ extends to an allowed eight-set $B=A\cup\{b\}$, which is
maximal, with $\sigma(B)=0$ and hence $\sigma(U)=b-a_0$. A representation of
$b-a_0$ by at most two terms of $2B$ is unavailable in $U$ only if it uses
$b$, which $U$ does not contain, or $a_0$ twice, so it is one of
\[
b,\qquad 2b,\qquad b+c\ \text{ with } c\in A,\qquad 2a_0 ,
\]
and these force $a_0=0$, $a_0=-b$, $c=-a_0$ and $b=0$ respectively. Each
contradicts the support conditions.
\end{proof}

\begin{proof}[Proof of Theorem~B for $m=3$]
The target lengths are $3t+7$. Degree zero is $\D(C_3^3)=7$;
degree one is Lemma~\ref{lem:ternary}(1).
 
For degree two, consider thirteen terms.
A zero set of at most three positions leaves at least ten terms, so
Lemma~\ref{lem:operations}(2) applies to the affine identity on the
remainder. Otherwise Lemma~\ref{lem:longternary} gives a zero set of at
least eleven positions. It properly contains a ten-term subsequence,
and Lemma~\ref{lem:operations}(3) raises that subsequence's affine
identity to degree two.
 
For degree three, a short zero set among sixteen terms again reduces to
the preceding degree. If there is none, the support bound and the
multiplicity bound force eight support points, each appearing twice.
Their sum is zero by Lemma~\ref{lem:ternary}(2), so the whole sequence
is a zero set. A thirteen-term subsequence properly inside it has a
degree-two identity, which Lemma~\ref{lem:operations}(3) raises to
degree three.
 
The support and multiplicity bounds also show that every seventeen
terms contain a short zero set: $\eta(C_3^3)\le17$ (in fact
$\eta(C_3^3)=17$, see \cite{ShortII}).
For $t\ge4$, the target length is at least nineteen, so the short-zero
recurrence proves all further degrees. Proposition~\ref{prop:sharp}
gives sharpness under its characteristic hypothesis.
\end{proof}

\subsection{The case \texorpdfstring{$C_4^3$}{C4 cubed}}
 
\begin{lemma}[Quaternary finite input]\label{lem:quaternary}
Every fourteen-term sequence over $C_4^3$ has an integral $\It_1$-identity,
and every eighteen-term sequence over $C_4^3$ has an integral
$\It_2$-identity.
\end{lemma}

The proof is computer-assisted; its reductions are as follows.
Since $\D(C_4^3)=10$, a zero set of at most four positions among
fourteen terms leaves a disjoint zero set. The two sets give an affine
identity by Lemma~\ref{lem:operations}(1). Each remaining fourteen-term
case receives either a disjoint-pair witness or an explicit integral
affine identity.
 
For eighteen terms, a short zero set reduces to the fourteen-term
claim by Lemma~\ref{lem:operations}(2).
Three disjoint zero sets also suffice, by part (1).
A zero set of at least fifteen positions gives a third route: it
properly contains fourteen positions, whose affine identity lifts by
part (3). The remaining cases are supplied with explicit integral
quadratic identities.
 
Appendix~\ref{app:computations} gives the complete finite domains,
the coverage arguments for two normalisations, and the counts of the positive
witnesses. The checks verify the group sums and every required
incidence moment over the integers. The enumeration includes
multiplicities and accounts for sequences in proper subgroups;
all transports preserve indexed occurrences.
 
The following example shows what the finite search can find beyond
these length-based reductions. Take
\[
U=(100)^3(010)^3(110)^3(001)^3(101)^3(131)^3
\]
in $(\Z/4\Z)^3$, with positions numbered $1,\dots,18$ in this order.
 
Inspecting the first, third and second coordinates shows that a
nonempty zero set of at most four terms would require four equal terms,
which are unavailable. Also $\sigma(U)=(0,3,1)$.
A zero set of at least fifteen terms would leave at most three terms
with this sum. Their first coordinates force them to be $010$ or
$001$, requiring three copies of $010$ and one of $001$, a contradiction.
Three disjoint zero sets would have a zero union of at least fifteen
terms, so they are impossible too.
 
There is nevertheless a useful five-term zero set
\[
B=\{9,11,13,14,16\}.
\]
Deleting it leaves only thirteen terms, one fewer than the length at
which an affine identity is guaranteed for every sequence.
For this particular remainder, such an identity can be written down.
The following five sets are zero sets disjoint from $B$:
\[
\begin{aligned}
A&=\{1,2,4,5,7,8\},\\
C&=\{3,8,10,12,15,18\},\\
D&=\{6,7,12,15,17,18\},\\
C'&=\{4,5,6,7,8,10,12,15,18\},\\
D'&=\{1,2,3,7,8,12,15,17,18\}.
\end{aligned}
\]
Position by position, they satisfy the integer identity
\[
\ind_A+\ind_C+\ind_D=\ind_{C'}+\ind_{D'}.
\]
The coefficients $1,1,1,-1,-1$ add to one and cancel through every
position, so they give an integral $\It_1$-identity on the remainder.
Lemma~\ref{lem:operations}(2) now produces an $\It_2$-identity on $U$:
keep these five weights, give their negatives to the unions of the
five sets with $B$, and give weight $1$ to $B$.
 
Thus the eleven-set quadratic certificate is a lower-degree identity
lifted by one of the elementary constructions. The useful extra
information is the affine identity on this particular thirteen-term
remainder; its existence does not follow from the universal
fourteen-term threshold. The certificate itself, in the form in which it
is stored and checked, is reproduced in the supplement.
 
\begin{proof}[Proof of Theorem~B for $m=4$]
The target lengths are $4t+10$.
Degree zero is $\D(C_4^3)=10$, and degrees one and two are
Lemma~\ref{lem:quaternary}.
The bound $\eta(C_4^3)\le22$ follows from the classical result
$\eta(C_4^3)=22$ of Harborth~\cite{Harborth}, see
\cite[Theorem~2.2(2)]{ShortII}.
For $t\ge3$, the target length is at least $22$; deleting a short
zero set leaves at least $4(t-1)+10$ terms.
Lemma~\ref{lem:operations}(2) therefore proves every further degree.
Proposition~\ref{prop:sharp} gives sharpness under its characteristic
hypothesis.
\end{proof}

\section{Davenport constants}\label{sec:davenport}

\begin{proof}[Proof of Corollary~C]
Theorems A and B give
\begin{equation}\label{eq:generalH}
\D(C_m^3\oplus H)\le m(\D(H)-1)+3m-2
\end{equation}
for every finite abelian $p$-group $H$.
Take $H=C_{p^a}\oplus C_{p^b}\oplus C_{p^c}$, allowing trivial
factors. Olson's theorem gives
\[
\D(H)-1=p^a+p^b+p^c-3.
\]
Since $p\nmid m$, the Chinese remainder theorem identifies
$C_m^3\oplus H$ with
$G=C_{mp^a}\oplus C_{mp^b}\oplus C_{mp^c}$.
Its invariant factors are in divisibility order, so substitution yields
\[
\D(G)\le m(p^a+p^b+p^c)-2=\DD(G).
\]
The basis lower bound gives the reverse inequality.
\end{proof}
 
The general bound~\eqref{eq:generalH} holds for kernels of every rank
and for every prime. Coprimality and the rank bound are used in the
comparison with $\DD$ above.
 
For example, take $m=3$ and $H=C_1\oplus C_5\oplus C_{25}$.
Here $\D(H)-1=28$, so the quotient identity is used in degree $28$,
at length $3(28+3)-2=91$. Transfer gives
\[
\D(C_3\oplus C_{15}\oplus C_{75})=91.
\]
No enumeration at length $91$ is needed: the constructions have already
extended the finite base cases to degree $28$. The basis sequence of
length $2+14+74=90$ shows why the resulting value is sharp.

\section{Verification, reproducibility and dependencies}\label{sec:data}

The proofs of Lemma~\ref{lem:ternary}(1) and
Lemma~\ref{lem:quaternary} are computer assisted. The support statement
in Lemma~\ref{lem:ternary}(2) follows from affine-plane counting,
and the proof of Lemma~\ref{lem:longternary} uses it.
The later incidence degrees are deductions from the finite base
identities by the three constructions. Table~\ref{tab:dependencies}
records these dependencies. Olson's small $p$-group values for the Davenport constant 
$\D(C_2^3)=4$, $\D(C_3^3)=7$ and $\D(C_4^3)=10$ also supply degree zero
and the short-zero reductions in the relevant rows.
 
For each computed case, the verification checks a positive witness:
disjoint zero sets, a suitable long zero set, or an integer incidence
identity. An incorrect negative classification, such as failing to
detect a disjoint pair, cannot by itself certify a case; a valid
positive witness must still pass the check.
Appendix~\ref{app:computations}
specifies the full finite domains and justifies the reductions used
to enumerate them.
The certificate verification uses exact arithmetic and does not rely
on the optimisation software used to discover the weights.
All certificates and the instructions for reproducing their results are
collected in a separate Supplement available on our  Github repository \cite{LeanRepo}.

The definitions, the lemmas, Proposition~\ref{prop:sharp},
Theorems~A and~B and Corollary~C have also been formalised
in Lean~4 with Mathlib~\cite{LeanRepo}, with proofs that
follow the ones given here.
In developing the formalisation, we made use of
Scott Armstrong's \emph{LeanAutoformalizationSkills}
collection~\cite{ArmstrongLeanSkills}.
For Lemma~\ref{lem:ternary}(1) and
Lemma~\ref{lem:quaternary}, the formalisation
performs its own enumerations; their completeness
is proved in Lean, together with the soundness
of the certificate checks (see the supplement).

\begin{table}[H]
\centering\small
\caption{Dependency map. ``Ordinary'' refers to the proof and its stated
inputs; ``computer-assisted'' means an exhaustive enumeration with
positive witnesses checked by an independent program; ``deduction''
means an ordinary argument some of whose inputs are computer-assisted
statements.}\label{tab:dependencies}
\begin{tabular}{>{\raggedright\arraybackslash}p{3.7cm}>{\raggedright\arraybackslash}p{2.5cm}>{\raggedright\arraybackslash}p{8.3cm}}
\toprule
Statement & Nature & Essential inputs\\
\midrule
Lemma~\ref{lem:duality} (duality) & ordinary & Boolean multilinearisation; linear algebra over a field\\
Lemma~\ref{lem:operations} (three constructions) & ordinary & Lemma~\ref{lem:duality}(1); inclusion--exclusion; two degree-cancellation arguments\\
Lemma~\ref{lem:indicator} (indicator) & ordinary & truncated group algebra; nilpotence of the augmentation ideal\\
Theorem~A & ordinary & evaluation identity; Lemma~\ref{lem:indicator}; Olson's value of $\D(H)$\\
Proposition~\ref{prop:sharp} (sharpness) & ordinary & integer triangular factorisation of the binomial-value matrix\\
Theorem~B, $m=2$ & ordinary & Fano identity; $\eta(C_2^3)\le8$; Lemma~\ref{lem:operations}; Proposition~\ref{prop:sharp}\\
Lemma~\ref{lem:ternary}(1) & computer-assisted & complete ten-term coverage (Appendix~\ref{app:computations}); positive affine or disjoint-pair certificates and their transports\\
Lemma~\ref{lem:ternary}(2) & ordinary & cap--plane incidence counts; linear functionals\\
Lemma~\ref{lem:longternary} & ordinary & maximal allowed supports; Lemma~\ref{lem:ternary}(2); multiplicity checks\\
Theorem~B, $m=3$ & deduction using finite input & Lemmas \ref{lem:ternary}, \ref{lem:longternary}, \ref{lem:operations}; $\D(C_3^3)=7$; $\eta(C_3^3)\le17$; Proposition~\ref{prop:sharp}\\
Lemma~\ref{lem:quaternary} & computer-assisted & two quaternary enumerations (Table~\ref{tab:quaternary}); exact certificates; the proper-subgroup bound used by the support-basis enumeration\\
Theorem~B, $m=4$ & deduction using finite input & Lemmas \ref{lem:quaternary}, \ref{lem:operations}; $\D(C_4^3)=10$; $\eta(C_4^3)\le22$; Proposition~\ref{prop:sharp}\\
Corollary~C & deduction using finite input & Theorems A and B; Olson's theorem; Chinese remainder theorem; basis lower bound\\
\bottomrule
\end{tabular}
\end{table}

\appendix
\section{Computational coverage and certificates}\label{app:computations}

This appendix collects the parts of the proofs of
Lemma~\ref{lem:ternary}(1) and Lemma~\ref{lem:quaternary} that justify
the finite computations: the exact domains that were enumerated and the
arguments showing that these domains cover every required sequence.
Nothing is inferred from samples: each enumeration establishes a
universal statement about all sequences of the stated length.
The positive witnesses themselves, their statistics, the stored form of
the eighteen-term example of Section~\ref{sec:bases} and the
instructions for an independent replay are collected in the
supplement.
 
\subsection{Ternary coverage}
 
A short-zero-free sequence over $C_3^3$ has all multiplicities at most two,
and its support $A$ contains no zero, no pair of opposite points and no three
distinct points with sum zero; equivalently, $A\cup\{0\}$ is an affine cap in
$\F_3^3$, a connection developed in~\cite{EdelEtAl}. A support with at least
five points spans $\F_3^3$, since a plane contains only four pairs of opposite
nonzero points. Mapping three independent support points to the standard
basis leaves ten pairs of opposite points, each contributing none or exactly
one of its two points, and the complete enumeration of these $3^{10}=59\,049$
choices produces the support counts of Table~\ref{tab:ternary}; no support of
size nine or more occurs. All $296$ seven-point supports have nonzero sum and
all $18$ eight-point supports have sum zero. A separate enumeration by direct
basis completion confirms these two facts.
These are the three facts of Lemma~\ref{lem:ternary}(2); they are
proved in Section~\ref{sec:bases} by counting affine planes, and the
enumeration is a redundant check of them.
 
For Lemma~\ref{lem:ternary}(1), a zero set of size at most three leaves seven terms and hence a
disjoint zero set, since $\D(C_3^3)=7$. It remains to treat the
short-zero-free multiplicity vectors of total ten; their exact counts are
given in Table~\ref{tab:ternary}. Most of them contain two disjoint zero
sets. The remaining $5\,006$ cases form $43$ orbits under the linear group,
and each orbit representative receives an explicit integral identity with
$5$, $7$ or $9$ zero sets and coefficients $\pm1$. An independent checker
regenerates the complete case space and verifies every zero set and every
moment of every identity, and a second one checks all transport matrices and
index permutations between representatives and cases. No affine translations
are used anywhere.
 
\begin{table}[ht]
\centering
\caption{The ternary ten-term computation, after normalising a support basis.
Support counts are not orbit counts.}\label{tab:ternary}
\begin{tabular}{rrrrr}
\toprule
Support size & Supports & Cases & Disjoint pair & Other identities\\
\midrule
5 &108&108&81&27\\
6 &296&4,440&2,769&1,671\\
7 &296&10,360&7,052&3,308\\
8 &18&504&504&0\\
\midrule
Total &&15,412&10,406&5,006\\
\bottomrule
\end{tabular}
\end{table}
 
\subsection{Quaternary coverage by reduction modulo two}
 
Represent $C_4^3$ as $(\Z/4\Z)^3$ and reduce modulo two. In a short-zero-free
sequence, the fibre over $0\in\F_2^3$ contains at most three terms, which
must be distinct and independent in $C_2^3$. A nonzero fibre contains at most
four distinct elements, since its eight elements form four pairs of opposite
elements, and at most one of these elements can repeat, because two distinct
repeated elements $g,h$ in the same fibre satisfy $2g+2h=0$. Each
multiplicity is at most three, so a fibre carries at most six terms.
 
We enumerate every projected multiplicity vector $(\mu_v)_{v\in\F_2^3}$ of the
prescribed total with $\mu_0\le3$ and $\mu_v\le6$ otherwise, up to the $168$
elements of $\GL_3(\F_2)$. Every such linear map lifts to an automorphism of
$(\Z/4\Z)^3$, because a binary lift has odd determinant. Write a lifted
element as $\widetilde v+2w$ with $v,w\in\F_2^3$, where $\widetilde v$ is the
binary representative. The automorphisms that reduce to the identity
modulo two, and their action on high parts, are
\[
\{\,I+2N:N\in M_3(\F_2)\,\},
\qquad w\longmapsto w+Nv,
\]
so choosing one term in each of the first independent occupied nonzero fibres,
a suitable $N$ makes each chosen high part $w$ zero. Since zero is the least
high part, requiring the first high part in each of these fibres to vanish,
after sorting high parts within fibres, retains a representative of every
orbit. This enumerator does not assume that the projection spans, so it covers
sequences contained in proper subgroups as well.
 
Backtracking then generates all lifts in nondecreasing order within each
fibre. It maintains the sets $R_1,R_2,R_3$ of sums of one, two and three
previously placed terms with distinct indices. Appending $x$ is legal exactly
when $x\ne0$ and $-x\notin R_1\cup R_2\cup R_3$, and the update is
\[
R'_1=R_1\cup\{x\},\quad
R'_2=R_2\cup(R_1+x),\quad
R'_3=R_3\cup(R_2+x),
\]
with the old sets on the right. This is an exact test for short-zero-freeness
that keeps track of repetitions.
 
\subsection{Independent support-basis enumeration}
 
A second enumeration normalises a support basis instead. An exhaustive
computation in the subgroup $C_4^2\oplus C_2$ shows that short-zero-free
sequences there have length at most thirteen. Every proper subgroup of $C_4^3$
lies in an index-two subgroup, all index-two subgroups are linearly
equivalent, and therefore every short-zero-free sequence of length at least
fourteen spans modulo two; three of its support vectors then have odd
determinant and form a basis of $(\Z/4\Z)^3$. We generate all short-zero-free
multisets containing the standard basis and canonicalise each by the least
sorted image under the maps that send an ordered support basis to the
standard basis; this is a complete invariant for linear equivalence. The two
enumerations produce identical canonical orbit lists.
 
\begin{table}[ht]
\centering\small
\caption{Quaternary finite coverage. Generated cases may repeat orbits; orbit
counts refer to the full linear group.}\label{tab:quaternary}
\begin{tabular}{lrr}
\toprule
Quantity & Length 14 & Length 18\\
\midrule
Generated by projection and lifts &308,051&18,224\\
Generated with a fixed support basis &24,445,377&1,267,272\\
Distinct linear orbits &106,655&4,941\\
Disjoint-zero witnesses &284,994&13,549\\
Long-zero witnesses (after disjoint cases)&---&4,171\\
Direct integer certificates &23,057&504\\
Orbits receiving direct certificates&6,963&64\\
\bottomrule
\end{tabular}
\end{table}

\end{document}